\documentclass[reqno]{amsart}%
\usepackage{amsmath}
\usepackage{amsfonts}
\usepackage{CJK}
\usepackage{amssymb}
\usepackage{graphicx}%
\usepackage{hyperref}
\usepackage{amscd}
\usepackage{mathrsfs}
\usepackage{bbm}

\newtheorem{thm}{Theorem}[section]

\theoremstyle{definition}
\newtheorem{defn}{Definition}[section]
\theoremstyle{Conjecture}

\theoremstyle{remark}
\newtheorem{rem}{Remark}[section]
\theoremstyle{Example}

\newcommand{\be}{\begin{equation}}
\newcommand{\ee}{\end{equation}}
\newcommand{\bea}{\begin{eqnarray}}
\newcommand{\eea}{\end{eqnarray}}
\newcommand{\ben}{\begin{eqnarray*}}
\newcommand{\een}{\end{eqnarray*}}
\newcommand{\bet}{\begin{equation}
\begin{split}}
\newcommand{\eet}{\end{split}
\end{equation}}

\begin{document}
\title[Analytic adjoint ideal sheaves]
{Analytic adjoint ideal sheaves associated to plurisubharmonic functions}

\author{Qi'an Guan}
\address{Qi'an Guan: School of Mathematical Sciences, and Beijing International Center for Mathematical Research,
Peking University, Beijing, 100871, China.}
\email{guanqian@amss.ac.cn}
\author{Zhenqian Li}
\address{Zhenqian Li: School of Mathematical Sciences, Peking University, Beijing, 100871, China.}
\email{lizhenqian@amss.ac.cn}

\thanks{The first author was partially supported by NSFC-11522101 and NSFC-11431013.}

\date{\today}

\begin{abstract}
In this article, we will present that the analytic adjoint ideal sheaves associated to plurisubharmonic functions are not coherent.
\end{abstract}

\maketitle

\section{Introduction}\label{sec:introduction}

The adjoint ideal sheaf on a smooth complex algebraic variety $X$ is a variant of the multiplier ideal sheaf in algebraic geometry (see \cite{Gue12, La04, Taka10} for more details).

In \cite{Gue12}, Guenancia gave an analytic definition of an adjoint ideal sheaf associated to a quasi-plurisubharmonic function $\varphi$ along a simple normal crossing (SNC) divisor $D=\sum D_i$ and established the compatibility with the algebraic adjoint ideal whenever $\varphi$ has analytic singularities.

Let $X$ be a complex manifold, $D=\sum D_i$ an SNC divisor and $\varphi$ a quasi-plurisubharmonic function on $X$. Let $Adj^{\alpha}_{D,*}(\varphi)\subset\mathcal{O}_X$ be the ideal sheaf of germs of holomorphic functions $f\in\mathcal{O}_{X,x}$ such that
$$|f|^2\prod_{k=1}^{p}\frac{1}{|h_k|^2(-\log|h_k|)^{\alpha}}e^{-\varphi}$$
is integrable with respect to the Lebesgue measure in some local coordinates near $x$, where $h=h_1\cdots h_p$ is the minimal defining function of $D$ near $x$ and $\alpha>1$.

In \cite{Gue12} (see also \cite{Kim15}), Guenancia gave the following analytic definition of adjoint ideal sheaf, which generalized the algebraic adjoint ideal sheaf (Proposition 2.11 in \cite{Gue12}; see also Proposition 5.1 in \cite{Kim15}).

\begin{defn} (\cite{Gue12, Kim15}).
The ideal sheaf $Adj^{\alpha}_{D}(\varphi):=\cup_{\varepsilon>0}Adj^{\alpha}_{D,*}((1+\varepsilon)\varphi)$ is called the \emph{analytic adjoint ideal sheaf} associated to $\varphi$ along $D$.
\end{defn}

Note that in \cite{Gue12} Guenancia used $\alpha=2$ in the definition, and later Kim in \cite{Kim15} extended the definition to $\alpha>1$ case. When $e^{\varphi}$ is locally H\"older continuous, Guenancia established the coherence of $Adj^{\alpha}_{D}$ for smooth divisor $D$ with $\varphi|_D\not\equiv-\infty$ (see Corollary 2.19 in \cite{Gue12}).

As mentioned by Guenancia and Kim, it is natural to ask the following\\
\\
\textbf{Question 1.1.} (\cite{Gue12, Kim15}). For $\alpha>1$ and $\varphi$ a general quasi-plurisubharmonic function on $X$, is the analytic adjoint ideal sheaf $Adj^{\alpha}_{D}(\varphi)$ coherent?\\

In this article, we will present the following negative answer to Question 1.1.

\begin{thm} \label{main}
There exists a plurisubharmonic function $\varphi$ and a smooth divisor $D$ with $\varphi|_D\not\equiv-\infty$ near the origin $o\in\mathbb{C}^n\ (n\geq3)$ such that for any $\alpha>1$, the analytic adjoint ideal sheaf $Adj^{\alpha}_{D}(\varphi)$ is not coherent near $o$.
\end{thm}

Specifically, we will construct $\varphi$ and $D$ with $\varphi|_D\not\equiv-\infty$ near $o$ such that the zero set of $Adj^{\alpha}_{D}(\varphi)$ is not an analytic set near $o$.

\section{Proof of main results}

We are now in a position to prove Theorem \ref{main}.\\
\\
\textbf{\emph{Proof of Theorem} \ref{main}.} Let $D=\{z_1=0\}$ and $$\varphi(z)=\max\{\sum\limits_{k=2}^{\infty}\alpha_k\log(|z_1|+|z_2-\frac{1}{k}|^{\beta_k}), \lambda\log|z_3|\},$$
where $\alpha_k=\frac{1}{2^{k!}},\ \beta_k=3\cdot2^{k!}$ and $\lambda>6$. Then $D$ is smooth and $\varphi(z)$ is a plurisubharmonic function near $o\in\mathbb{C}^n\ (n\geq3)$. Without loss of generality, we assume $n=3$ and $U$ contained in the unit polydisk $\Delta^3$ is a neighborhood of $o$ such that $\log(|z_1|+|z_2-\frac{1}{k}|^{\beta_k})<0$ for any $k\geq2$ on $U$.\\

Step 1. \emph{Nonintegrability of $\frac{1}{|z_1|^2(-\log|z_1|)^{\alpha}}e^{-\varphi}$ near $(0,\frac{1}{k},0)$ for any $(0,\frac{1}{k},0)\in U$.}

Since
\begin{equation*}
\begin{split}
\varphi(z)&=\max\{\sum\limits_{k=2}^{\infty}\alpha_k\log(|z_1|+|z_2-\frac{1}{k}|^{\beta_k}), \lambda\log|z_3|\}\\
&\leq\max\{\alpha_k\log(|z_1|+|z_2-\frac{1}{k}|^{\beta_k}), \lambda\log|z_3|\}\\
&\leq\log\big((|z_1|+|z_2-\frac{1}{k}|^{\beta_k})^{\alpha_k}+|z_3|^{\lambda}\big)
\end{split}
\end{equation*}
near $o$, replacing $z_2-\frac{1}{k}$ by $z_2$ near $(0,\frac{1}{k},0)$, for sufficiently small polydisk $\Delta_r^3$ we have
\begin{equation*}
\begin{split}
&\int_{\Delta_r^3}\frac{1}{|z_1|^2(-\log|z_1|)^{\alpha}}e^{-\varphi(z_1,z_2+\frac{1}{k},z_3)}dV_{3}\\
\geq&\int_{\Delta_r^3}\frac{dV_{3}}{\big((|z_1|+|z_2|^{\beta_k})^{\alpha_k}+|z_3|^{\lambda}\big)
 |z_1|^2(-\log|z_1|)^{\alpha}}\\
=&\int_{\Delta_r^*}\frac{dV_1}{|z_1|^2(-\log|z_1|)^{\alpha}}
 \int_{\Delta_r^2}\frac{|z_1|^{\frac{2}{\beta_k}+\frac{2\alpha_k}{\lambda}}}{|z_1|^{\alpha_k}}
 \frac{(\frac{\sqrt{-1}}{2})^2d(\frac{z_2}{|z_1|^{\frac{1}{\beta_k}}})\wedge d(\frac{\bar z_2}{|z_1|^{\frac{1}{\beta_k}}})\wedge d(\frac{z_3}{|z_1|^{\frac{\alpha_k}{\lambda}}})\wedge d(\frac{\bar z_3}{|z_1|^{\frac{\alpha_k}{\lambda}}})}
 {(1+\frac{|z_2|^{\beta_k}}{|z_1|})^{\alpha_k}+\frac{|z_3|^{\lambda}}{|z_1|^{\alpha_k}}}\\
\geq&\int_{\Delta_r^*}\frac{dV_1}{|z_1|^{2+(\alpha_k-\frac{2}{\beta_k}-\frac{2\alpha_k}
 {\lambda})}(-\log|z_1|)^{\alpha}}
 \int_{\Delta}\frac{(\frac{\sqrt{-1}}{2})^2dw_2\wedge d\bar w_2\wedge dw_3\wedge d\bar w_3}{(1+|w_2|^{\beta_k})^{\alpha_k}+|w_3|^{\lambda}}\\
\geq&\ C\cdot\int_{\Delta_r^*}\frac{dV_1}{|z_1|^{2+(\alpha_k-\frac{2}{\beta_k}-\frac{2\alpha_k}{\lambda})}
 (-\log|z_1|)^{\alpha}}=+\infty,
\end{split}
\end{equation*}
where $C>0$ is some constant and the last equality is from $\alpha_k-\frac{2}{\beta_k}-\frac{2\alpha_k}{\lambda}>0$. It follows that $\frac{1}{|z_1|^2(-\log|z_1|)^{\alpha}}e^{-\varphi}$ is not locally integrable near $(0,\frac{1}{k},0)$ for any $(0,\frac{1}{k},0)\in U$.\\

Step 2. \emph{Integrability of $\frac{1}{|z_1|^2(-\log|z_1|)^{\alpha}}e^{-(1+\varepsilon)\varphi}$ near $(0,z_2,0)\in U$ with $z_2\not=\frac{1}{k},0$.}

Take $0<\varepsilon_0<1$ such that $\alpha-\varepsilon_0>1$ and $|z_2-\frac{1}{k}|>\varepsilon_0$ for any $k$. Then, for any $N\geq1$ with $\varepsilon_0^{\beta_{N+1}}<|z_1|\leq\varepsilon_0^{\beta_{N}}$, we obtain
\begin{equation*}
\begin{split}
\varphi(z)&\geq\sum\limits_{k=2}^{\infty}\alpha_k\log(|z_1|+|z_2-\frac{1}{k}|^{\beta_k})\geq\sum\limits_{k=2}^{\infty}
 \alpha_k\log(|z_1|+\varepsilon_0^{\beta_k})\\
&\geq\sum_{k\leq N}\alpha_k\log(|z_1|+\varepsilon_0^{\beta_k})
 +\sum_{k\geq N+1}\alpha_k\log(|z_1|+\varepsilon_0^{\beta_k})\\
&\geq\sum_{k\leq N}\alpha_k\beta_k\log\varepsilon_0+\sum_{k\geq N+1}\alpha_k\log|z_1|
 \geq3N\log\varepsilon_0+\sum_{k\geq N+1}\alpha_k\beta_{N+1}\log\varepsilon_0\\
&\geq3N\log\varepsilon_0+2\alpha_{N+1}\beta_{N+1}\log\varepsilon_0=3(N+2)\log\varepsilon_0
\end{split}
\end{equation*}
and
\begin{equation*}
\begin{split}
&\log(-\log|z_1|)\geq\log(-\beta_{N}\log\varepsilon_0)=N!\log2+\log3+\log(-\log\varepsilon_0).\qquad (*)
\end{split}
\end{equation*}

Set $C_N:=N!\log2+\log3+\log(-\log\varepsilon_0)$. Since $$0<\frac{3(N+2)\log\varepsilon_0}{-C_N}\to0\ (N\to\infty),$$
it follows from $(*)$ that for large N, we have 
$$\frac{3(N+2)\log\varepsilon_0}{-\log(-\log|z_1|)}\leq\frac{3(N+2)\log\varepsilon_0}{-C_N}
\leq\frac{\varepsilon_0}{1+\varepsilon}.$$
Hence, for $\varepsilon_0^{\beta_{N+1}}<|z_1|\leq\varepsilon_0^{\beta_{N}}$, we obtain $$(1+\varepsilon)\varphi\geq3(1+\varepsilon)(N+2)\log\varepsilon_0\geq-\varepsilon_0\log(-\log|z_1|)$$ for large enough $N$. Thus, when $|z_1|$ is small enough, we have
$$\frac{1}{|z_1|^2(-\log|z_1|)^{\alpha}}e^{-(1+\varepsilon)\varphi}
\leq\frac{1}{|z_1|^2(-\log|z_1|)^{\alpha-\varepsilon_0}},$$
which is locally integrable near $(0,z_2,0)\in U$ with $z_2\not=\frac{1}{k},0$ by $\alpha-\varepsilon_0>1$.\\

Step 3. \emph{Incoherence of the analytic adjoint ideal sheaf $Adj^{\alpha}_{D}(\varphi)$ near $o$.}

Suppose that $Adj^{\alpha}_{D}(\varphi)$ is a coherent ideal sheaf on $U$. Then, the zero set $N(Adj^{\alpha}_{D}(\varphi))$ of $Adj^{\alpha}_{D}(\varphi)$ is an analytic set in $U$. However, it follows from Step 1 and Step 2 that on $U$,
$$N(Adj^{\alpha}_{D}(\varphi))=\{(0,\frac{1}{k},0)\}\cup\{o\},$$
which is not analytic at $o$, contradicting to the assumption.
\hfill $\Box$

\begin{rem}
In Theorem \ref{main}, if the condition $\varphi|_D\not\equiv-\infty$ is not required, then the result also holds for $n=2$ by a slight modification. Indeed, we can take $D=\{z_1=0\}$ and $\varphi(z_1,z_2)=\sum\limits_{k=2}^{\infty}\alpha_k\log(|z_1|+|z_2-\frac{1}{k}|^{\beta_k})$ with $\alpha_k=\frac{1}{2^{k!}},\ \beta_k=3\cdot2^{k!}$.
\end{rem}


\begin{thebibliography}{123}
\bibitem{Gue12} H. Guenancia, \emph{Toric plurisubharmonic functions and analytic adjoint ideal sheaves}, Math. Z. 271 (2012), 1011--1035.
\bibitem{Kim15} D. Kim, \emph{Themes on Non-analytic Singularities of Plurisubharmonic Functions}, Complex Analysis and Geometry, Springer Proceedings in Mathematics $\&$ Statistics, vol. 144, Springer, Tokyo, 2015, pp. 197--206.
\bibitem{La04} R. Lazarsfeld, \emph{Positivity in Algebraic Geometry II}, Springer, New York, 2004.
\bibitem{Taka10} S. Takagi, \emph{Adjoint ideals along closed subvarieties of higher codimension}, J. Reine Angew. Math. 641 (2010), 145--162.
\end{thebibliography}
\end{document}